\documentclass[12pt, reqno]{amsart}
\usepackage{hyperref}
\usepackage{amssymb,amsmath,graphicx,amsthm}
\def\simgeq{\underset\sim>}

\def\1#1{\overline{#1}}

\def\2#1{\widetilde{#1}}

\def\3#1{\widehat{#1}}

\def\4#1{\mathbb{#1}}

\def\5#1{\frak{#1}}

\def\6#1{{\mathcal{#1}}}

\def\C{{\4C}}

\newcommand{\Om}{\Omega}

\def\C{{\Bbb C}}

\def\di{\partial}
\def\dib{\bar\partial}

\def\Label#1{\label{#1}}

\def\simgeq{\underset\sim>}

\def\1#1{\overline{#1}}

\def\2#1{\widetilde{#1}}

\def\3#1{\widehat{#1}}

\def\4#1{\mathbb{#1}}

\def\5#1{\frak{#1}}

\def\6#1{{\mathcal{#1}}}

\def\C{{\4C}}

\numberwithin{equation}{section}

\theoremstyle{plain}

\newtheorem{theorem}{Theorem}[section]

\theoremstyle{definition}

\theoremstyle{remark}

\begin{document}
\begin{abstract}
We give a simple proof of the fact that an ``$f$-estimate" for the $\dib$-Neumann problem implies a lower bound on the geomatric type of the boundary along any complex one dimensional variety. The proof uses the existence of  peak functions which is in turn a consequence of  the $f$-estimate.
\end{abstract}
\title[Lower bound for the geometric type]{Lower bound for the geometric type from a generalized estimate in the $\dib$-Neumann problem - a new approach by peak functions}
  
  \author[T.V. Khanh]{Tran Vu Khanh}   
\maketitle
\section{Introduction}
In a series of seminal papers in the Annal of Mathematics \cite{Cat83,Cat87}, Catlin proved the equivalence of the finite type  of a boundary (cf. \cite{Dan82}) with the existence of a subelliptic estimate for  the $\dib$-Neumann problem by triangulating through the $t^\epsilon$-Property (see below)
\begin{enumerate}
  \item[(i)] finite type $m$ $\Rightarrow$ $t^\epsilon$-Property with $\epsilon=m^{-n^2m^{n^2}}$;
  \item[(ii)] $t^\epsilon$-Property $\Rightarrow$ $\epsilon$-subelliptic estimate;
  \item[(iii)] $\epsilon$-subelliptic estimate $\Rightarrow$  finite type $m$ for $m\le\frac{1}{\epsilon}$.
\end{enumerate}
Here, the $t^\epsilon$-Property  of a boundary $b\Om$ is a case of a more general ``$f$-Property"  defined  as follows. For a smooth monontonic increasing function $f :[1+\infty)\to[1,+\infty)$ with $f (t)\le t^{\frac{1}{2}}$,  the $f$-Property at $z_o$ consists in the existence of a neigborhood $U$ of $z_o$, of constant $C_1, C_2$ and of a family of functions $\{\phi_\delta\}$ such that
\begin{enumerate}
  \item [1)] $\phi_\delta$ are plurisubharmonic and $C^2$ on $U$ and $-1\le \phi_\delta \le0$;
  \item[2)] $\di\dib \phi_\delta\ge C_1f(\delta^{-1})^2Id$ and $|D\phi_\delta|\le C_2  \delta^{-1}$ for any  $z\in U\cap \{z\in \Om:-\delta<r(z)<0\}$, where $r$ is a defining function of $\Om$.
  \end{enumerate}
Also, the finite type $m$ is a case of a general ``$F$-type" along a one-dimensional complex analytic variety $Z$ consisting in the estimate
\begin{eqnarray}\Label{aa}
|r(z)|\le F(|z-z_o|) \quad z\in Z, z\to z_o.
\end{eqnarray}
The results in the step (ii) and (iii) have been generalized in \cite{KZ10, KZ12}. In particular, in \cite{KZ10}  it has been showed that the $f$-Property implies an $f$-estimate  for any $f$, and in \cite{KZ12}  that an  $f$-estimate with $\dfrac{f}{\log }\to \infty$ at $\infty$  implies that the type along complex analytic variety has a lower bound with the rate $G$ with 
\begin{eqnarray}\Label{G1}
G(\delta)=\left(\left(\frac{f}{\log}\right)^*(\delta^{-1})\right)^{-1},
\end{eqnarray}
where the superscript $^*$ denotes the inverse function. Combining the above results, one obtains
\begin{theorem}[Catlin \cite{Cat83, Cat87}; Khanh-Zampieri \cite{KZ10, KZ12}]\Label{B} Let $\Om$ be a pseudoconvex domain in $\C^n$ with $C^\infty$-smooth boundary and $z_o$ be a point in the boundary $b\Om$. Assume that the $f$-Property holds at $z_o$ for $\dfrac{f}{\log }\nearrow \infty$ as $t\to\infty$ and that $b\Om$ has type $\le F$ along an one-dimensional complex analytic variety at $z_o$.
Then, $F\simgeq G$.
\end{theorem}

 The purpose of this note is to give an immediate  proof of Theorem~\ref{B} which has also the advantage of requiring a  minimal smoothness of $b\Om$ and the expenses of a slightly stronger assumption about $f$.  In fact, we prove the following
\begin{theorem}\Label{t1}Let $\Om$ be a pseudoconvex domain of $\C^n$ with $C^2$-smooth boundary $b\Om$ and $z_o$ be a point in the boundary $b\Om$. Assume that $\Om$ has the $f$-Property at $z_o$ with $f$ satisfying $\displaystyle\int_t^\infty \dfrac{da}{a f(a)}<\infty$ for some $t>1$,  denote by $(g(t))^{-1}$ this finite integral, and set $G(\delta)=\left(g^*(\delta^{-1})\right)^{-1}$.  Then, if $b\Om$ has type $\le F$ along  one-dimensional complex analytic variety $Z$, then $F(\delta)\simgeq G(\alpha\delta)$ for a suitable $\alpha$ and for any $\delta$ small.

\end{theorem}
First, we remark that the $C^\infty$-smoothness of the boundary in the results of Catlin and Khanh-Zampieri (Theorem \ref{B}) is required since they apply the regularity of the $\dib$-Neumann problem. In Theorem~\ref{t1}, the  condition of smoothness is reduced because of the use of a plurisubharmonic peak function. However, in the construction of the family of the plurisubharmonic peak functions, we need a slightly stronger hypothesis in $f$, for example, if $f(t)=\log t\cdot\log^\epsilon(\log t)$ with $0<\epsilon\le 1$, then $f$ satisfies the hypothesis in Theorem~\ref{B} but does not in Theorem~\ref{t1}. Finally, we notice the equivalence of the statements in the two theorems in case $f(t)=\log^\beta t$ for $\beta>1$, or $f(t)=t^\epsilon$ for any $0<\epsilon\le \frac{1}{2}$.
\section{Proof of Theorem~\ref{t1}}
The proof of Theorem~\ref{t1} follows immediately from Theorem~\ref{pshpeak} and Theorem~\ref{lowerbound}  below.  In \cite{Kha13}, the author proves that there exists a family of plurisubharmonic functions with good estimates 

\begin{theorem}\Label{pshpeak} Under the assumptions of Theorem~\ref{t1},  for a fixed constant $0< \eta\le 1$, there are a neighborhood $V$ of $z_o$ and positive constants $c_1, c_2, c_3$ such that the following holds.  For any $w\in V\cap b\Om$ there is a plurisuhharmonic function $\psi_w$ on $V\cap\Om$ verifying
 \begin{enumerate}
   \item $|\psi_w(z)-\psi_w(z')|\le c_1|z-z'|^\eta$,
      \item $\psi_w(z)\le-G^\eta(c_2|z-w|)$,
   \item $\psi_{\pi(z)}(z)\ge-c_3\delta_{b\Om}(z)^{\eta}$,
 \end{enumerate}
for any $z$ and $z'$ in $V\cap\bar\Om$ (where $\delta_{b\Om}(z)$ and $\pi(z)$ denote the distance and projection of $z$ to the boundary, respectively).
\end{theorem}
Using the conclusion of Theorem \ref{pshpeak} just for $w\equiv z_0$, we get 
\begin{theorem}\Label{lowerbound}
Let $\Om$ be a $C^2$-smoothly pseudoconvex domain in $\C^n$ and $z_o$ be a boundary  point. Assume that there is a neighborhood $V$ of $z_o$ and a plurisubharmonic function $\psi$ on $V\cap \Om$ such that    
\begin{eqnarray}\Label{key} 
   -c_1|z-z_o|^\eta\le \psi(z)\le -G^\eta(c_2|z-z_o|),\quad z\in V\cap\Om,
\end{eqnarray}
for suitable $c_1, c_2>0$ and $\eta\in (0,1]$.  In this situation, if $b\Om$ has type $\le F$ along a one-dimensional complex analytic variety $Z$, then $F(\delta)\simgeq  G(\alpha\delta)$.
\end{theorem}
In the case $G(t)=t^m$, the result has been obtained by Fornaess and Sibony in \cite{FS89}.\\
{\it Proof of Theorem~\ref{lowerbound}.} Let $\Om$ be a domain in $\C^n$ and assume that there is a function $F$ and an $1$-dimensional complex analytic variety $Z$ passing through $z_o$ such that \eqref{aa} is satisfied for $z\in Z$. Then, in any neighborhood $U$ of $z_o$ there are  constants $c_3, c_4>0$ and a family $\{Z_\delta\}$ of $1$-dimensional complex manifolds $Z_\delta\subset U$ defined by $h_\delta:\overline\Delta\to U$ with $h_\delta(0)=z_o$ such that  
\begin{eqnarray}\Label{dd}
\delta=\sup_{t\in \overline\Delta}|h_\delta(t)-z_o|\ge |h_\delta'(0)| \ge c_3\delta
\end{eqnarray}
and 
\begin{eqnarray}\Label{cc}
\sup_{t\in \overline\Delta}|\delta_{b\Om}(h_\delta(t))|< c_4F(\delta),
\end{eqnarray}
where  $\Delta$ denotes the unit disc in $\C$.\\

Let $\nu$ be the outward normal to $b\Om$ at $z_o$. From \eqref{cc}, we have $h_\delta(t)-c_4F(\delta)\nu\in \Om\cap U$ for any $t\in\overline\Delta$. Applying the submean value inequality to the subharmonic function $\psi(h_\delta(t)-c_4F(\delta)\nu)$ on $\overline\Delta$ we get
\begin{eqnarray}\Label{bb}
\begin{split}
\psi(z_o-c_4F(\delta)\nu)\le& \frac{1}{2\pi}\int_{0}^{2\pi}\psi(h_\delta(e^{i\theta})-c_4F(\delta)\nu)d\theta.
\end{split}
\end{eqnarray}
Now, we use the first inequality in \eqref{key} for the left hand side term of \eqref{bb}
$$-\psi(z_o-c_4F(\delta)\nu)\le c_1c_4^\eta F^\eta(\delta)^\eta.$$
For the right hand side term of \eqref{bb}, we use the second inequality of \eqref{key}
\begin{eqnarray}\Label{ee}
\begin{split}
-\frac{1}{2\pi}\int_{0}^{2\pi}\psi(h_\delta(e^{i\theta})-c_4F(\delta)\nu)d\theta\ge &\frac{1}{2\pi}\int_{0}^{2\pi}G^\eta\left(c_2\left|h_\delta(e^{i\theta})-c_4F(\delta)\nu-z_o\right|\right)d\theta.
\end{split}
\end{eqnarray}
Using \eqref{dd} and the Jensen inequality for the increasing, convex function $G^\eta$, we get 
\begin{eqnarray}\Label{ff}
\begin{split}
G^\eta(c_2c_3\delta)\le G^\eta(c_2 |h_\delta'(0)|)
\le &G^\eta \Big(\frac{1}{2\pi}\int^{2\pi}_0c_2|h_\delta(e^{i\theta})-CF(\delta)\nu-z_o|d\theta\Big)\\
\le&\frac{1}{2\pi} \int^{2\pi}_0G^\eta\big(c_2|h_\delta(e^{i\theta})-CF(\delta)\nu-z_o|\big)d\theta .
\end{split}
\end{eqnarray}
Combining \eqref{bb}, \eqref{ee} and \eqref{ff}, we obtain 
$$F(\delta)\ge \beta G(\alpha\delta).$$
with $\beta=\left(c_1^{1/\eta}c_4\right)^{-1}$ and $\alpha=c_2c_3$. The proof of Theorem \eqref{t1} is complete.

$\hfill\Box$

\bibliographystyle{alpha}

\begin{thebibliography}{Koh79}

\bibitem[Cat83]{Cat83}
David Catlin.
\newblock Necessary conditions for subellipticity of the {$\bar \partial
  $}-{N}eumann problem.
\newblock {\em Ann. of Math. (2)}, 117(1):147--171, 1983.


\bibitem[Cat87]{Cat87}
D. Catlin.
\newblock Subelliptic estimates for the {$\overline\partial$}-{N}eumann problem
  on pseudoconvex domains.
\newblock {\em Ann. of Math. (2)}, 126(1):131--191, 1987.

\bibitem[D'A82]{Dan82}
John~P. D'Angelo.
\newblock Real hypersurfaces, orders of contact, and applications.
\newblock {\em Ann. of Math. (2)}, 115(3):615--637, 1982.

\bibitem[FS89]{FS89}
J~E. Fornaess and N. Sibony.
\newblock Construction of P.S.H. functions on weakly pseudoconvex domains.
\newblock{\em Duke Math. J.} 58 (3): 633--655, 1989. 

\bibitem[Kha13]{Kha13}
T.~V Khanh.
\newblock Boundary behavior of the Kobayashi metric near a point of infinite type.
\newblock 2013.

\bibitem[KZ10]{KZ10}
T.~V Khanh and G. Zampieri.
\newblock Regularity of the {$\overline\partial$}-{N}eumann problem at point of
  infinite type.
\newblock {\em J. Funct. Anal.}, 259(11):2760--2775, 2010.

\bibitem[KZ12]{KZ12}
T.~V Khanh and G. Zampieri.
\newblock Necessary geometric and analytic conditions for general estimates in
  the {$\bar{\partial}$}-{N}eumann problem.
\newblock {\em Invent. Math.}, 188(3):729--750, 2012.



\end{thebibliography}

\end{document}